\documentclass[11pt,a4paper,underruledhead]{paper}
\usepackage{enumerate,math}
\usepackage[mathscr]{euscript}
\usepackage{float,eepic}

\title[A note on orthogonality of subspaces in Euclidean geometry]{A note on orthogonality of subspaces in\\ Euclidean geometry}
\author{J. Konarzewski,  M. {\.Z}ynel}

\newfloat{rysunek*}{t}{flt}[section]
\floatname{rysunek*}{Figure}

\newenvironment{rysunek}[1]{
  \def\figlabel{#1}
  \begin{rysunek*}
    \begin{center}
      \small
}{%
    \caption{}
    \label{\figlabel}
    \end{center}
  \end{rysunek*}}

\newsavebox\tempboxa

\expandafter\def\csname floatc@plain\endcsname#1#2{%
  \let\captionsep\relax
  \setbox\tempboxa\hbox{#2}
  \ifdim\wd\tempboxa>0pt
    \def\captionsep{:\space\space}
  \fi
  \setbox\tempboxa\hbox{#1\captionsep #2}%
  \ifdim\wd\tempboxa>\hsize #1\captionsep #2\par
    \else\hbox to\hsize{\hfil\box\tempboxa\hfil}\fi}

\def\lines{{\cal L}}
\def\hipy{{\cal H}}

\def\perpx{\mathrel{\perp^{\mkern-4mu\ast}}}
\def\perpgo{\mathrel{\perp\mkern-16.5mu\raise0.2ex\hbox{$\circ$}}}
\def\perpg{\mathrel{\perp\mkern-19mu\raise-0.76ex\hbox{$=$}}}
\def\perppo{\mathrel{\perp\mkern-11mu\perp}}
\def\perpu#1#2#3{\perp^{#1}_{{#2},{#3}}}
\def\collin{\sim}

\def\oadjac{\mathrel{\raise-0.5ex\hbox{$\collin$}\mkern-14mu{\raise0.4ex\hbox{$\perp$}}}}

\newenvironment{ctext}{%
  \par
  \smallskip
  \centering
}{%
 \par
 \smallskip
 \csname @endpetrue\endcsname
}

\def\defiff{:\mkern-3mu\iff}

\begin{document}

\maketitle

\begin{abstract}
  We show that Euclidean geometry in suitably high dimension
  can be expressed as a theory of orthogonality of subspaces
  with fixed dimensions and fixed dimension of their meet.

  \begin{flushleft}
  \end{flushleft}  
\end{abstract}

\section{Introduction}

While the notion of orthogonality of lines in Euclidean geometry has well 
founded meaning (it is frequently used as a primitive notion, see \cite{lenz}),
orthogonality of subspaces can be defined in several different ways.
Two of them were shown in \cite{perpadja} to be sufficient in Euclidean geometry;
actually, each of these two considered on the universe of subspaces of fixed 
dimension can be used to reinterpret the underlying point-line affine space
and after that to define line orthogonality.
Thus the procedure of reinterpretation consists, in fact, in two steps
and in the second step one should define orthogonality of lines in terms of 
a given orthogonality of subspaces. 
In this note we show that such a definition is possible for each prescribed values
of dimensions of the considered subspaces (Theorem~\ref{thm:main}\eqref{main:2}).

The notion of orthogonality of subspaces is not a unique-meaning relation,
even if dimensions of the subspaces involved are fixed.
Therefore, we have to deal with a family of possible relations of 
orthogonality. 
And in this note we show that {\em each one} of these relations is sufficient to 
express the underlying geometry provided the latter has sufficiently high 
dimension (Theorem~\ref{thm:main}\eqref{main:1}).

So, finally, we prove that Euclidean geometry can be expressed in the language
with points, subspaces (of fixed dimensions), and orthogonality of subspaces.
It is a folklore that affine geometry can be expressed as a theory of 
point-$k$-subspace incidence. Euclidean geometry appears when we impose a
relation of orthogonality on that ``affine'' structure.

Our result does not solve the problem whether Euclidean geometry can be expressed
in the language with $k$-subspaces {\em as individuals} and some 
of the orthogonalities introduced above as a single primitive notion, in that way,
possibly, generalizing \cite{perpadja}. We conjecture that the answer is affirmative, 
but the question is addressed in other papers.

We close the paper with a list of some more interesting properties of 
the orthogonalities considered here. This list is not intended as a complete
axiom system, but we think that at least some of its items can be used to build
such a system characterizing orthogonality of subspaces.

\section{Results}

Let ${\goth M} = \struct{S,\lines,\perp}$ be an Euclidean space, where
${\goth A} := \struct{S,\lines}$ is an affine space with $\lines\subset 2^S$
and $\mathord{\perp}\subset \lines\times\lines$ is a line orthogonality
(cf. \cite{lenz}).
Up to an isomorphism ${\goth M}$ corresponds to $\struct{V, \lines_V, \perp_\xi}$ 
where $V$ is a vector space, $\lines_V$ is the set of translates of 
1-dimensional subspaces of $V$ and $\perp_\xi$ is the orthogonality determined
by a nondegenerate symmetric bilinear form $\xi$ on $V$ with no isotropic directions.
For each nonnegative integer $k$, 
$\hipy_k$ stands for the class of all $k$-dimensional subspaces of $\goth M$, 
and $\hipy$ stands for all subspaces of $\goth M$.
If $X_1,X_2\in\hipy$ we write 
$X_1 \sqcup X_2$ for the least subspace in $\hipy$
that contains
$X_1\cup X_2$ (i.e. the meet of all elements of $\hipy$ containing
$X_1\cup X_2$).
Note an evident fact that follows from elementary affine geometry.
\begin{fact}\label{fct:hipy<->lines}
  \begin{sentences}\itemsep-2pt
  \item
    The family $\hipy_k$ is definable in $\goth A$ for each nonnegative integer $k$.
  \item
    Let $k < \dim({\goth A})$. Then the family $\lines$ is definable in the incidence
    structure $\struct{S,\hipy_k}$.
    Consequently, $\goth A$ is definable in $\struct{S,\hipy_k}$.
  \end{sentences}
\end{fact}
Recall that $\goth M$ is definitionally equivalent to the structure 
$\struct{S,\lines,\perppo}$ (cf. e.g. an axiom system for $\perppo$ in \cite{kusak}, \cite{szczerba1}),
where $\mathord{\perppo}\subset S^2\times S^2$
is defined in $\goth M$ by the formula
\begin{equation}\label{def:perp:pt}
  a,b \perppo c,d \defiff \text{ there are } L_1,L_2\in\lines \text{ such that }
  a,b \in L_1 \perp L_2 \ni c,d.
\end{equation}
Given any two $X,Y\in\hipy$ we write 
\begin{equation}\label{def:perp}
  X\perp Y \defiff  a,b \perppo c,d \text{ for all } a,b\in X,\; c,d\in Y.
\end{equation}
Note that for $X,Y\in\lines$ the relation defined by \eqref{def:perp}
coincides with the orthogonality we have started from.
If $X\perp Y$ then $X\cap Y$ is at most a point;
we write
\begin{equation}\label{def:perpx}
  X\perpx Y \defiff X\perp Y \text{ and } X\cap Y\neq \emptyset.
\end{equation}

\begin{rysunek}{fig:perpg}
  \setlength{\unitlength}{0.00037in}
  \begin{picture}(8577,6039)(0,-10)
    \path(12,12)(12,4812)(12,4812)(3162,6012)
    \path(3162,1212)(8262,1212)
	    (5112,12)(12,12)
    \put(3194.268,1312.291){\arc{1301.019}{4.6628}{6.4380}}
    \put(3739.500,934.500){\arc{2357.340}{3.1861}{4.2003}}
    \put(3025.636,2084.727){\arc{2383.244}{0.8218}{1.9704}}
    \put(3482,1540){\circle*{42}}
    \put(2902,1427){\circle*{42}}
    \put(3192,1062){\circle*{42}}
    \put(3162,1212){\circle*{84}}
    \put(237,4362){$X_1$}
    \put(4812,237){$X_2$}
    \put(3237,4662){$Z_1$}
    \put(6612,1362){$Z_2$}
    \put(512,800){$X_1\cap X_2$}
    \put(3312,1260){$q$}
    \allinethickness{1.2pt}
    \path(3162,1212)(12,12)
    \path(3162,1212)(3162,6012)
    \path(3162,1212)(8262,1212)   
  \end{picture}
\end{rysunek}

Recall that for any two subspaces $X, V\in\hipy$ such that
$X\subset V$ and a point $q\in X$ there is the unique maximal $X'\in\hipy$
such that $q \in X'\perpx X$, \, $X'\subset V$,
and $X\sqcup X' = V$.
We call $X'$ an \emph{orthocomplement} of $X$ in $V$ through $q$.
If $X$ is a point then necessarily $X=\{q\}$ and $X' = V$.

Let us define now (cf. Figure~\ref{fig:perpg})
\begin{multline}\label{def:perpgo}
  X_1\perpgo X_2\;\; \defiff
    \text{ there is a point } q \in X_1\cap X_2 \text{ and } Z_1, Z_2\in\hipy \text{ such that } \\
      q\in Z_1, Z_2 \perpx X_1 \cap X_2, \;\; 
      Z_1 \perpx Z_2 \text{ and } 
      (X_1 \cap X_2) \sqcup Z_i = X_i \text{ for } i=1,2.
\end{multline}
It is seen that the relation $\perpgo$ is symmetric. It is also not too hard to note
that the following holds
\begin{equation}\label{wz:osieo}
  X_1 \perpgo X_2 \iff \text{there is } Z_i\in\hipy \text{ such that } 
  Z_i \perpx X_{3-i} \text{ and }  (X_1\cap X_2)\sqcup Z_i = X_i
\end{equation}
for both $i=1,2$.
Note that when $X_1\cap X_2$ is a point then 
$X_1\perpg X_2$ and $X_1\perpx X_2$ are equivalent.
Recall also a known formula
\begin{equation}\label{eq:perpxsup}
  q \in X_1,X_2 \perpx Y\ni q \implies  Y\perpx (X_1\sqcup X_2).
\end{equation}

The motivation for such general definition \eqref{def:perpgo} is reflection 
geometry (cf. \cite{bachman}, \cite{schroeder}). 
Denote by $\sigma_X$ the reflection in a subspace $X$, i.e. an involutory 
isometry that fixes $X$ pointwise; then
\begin{equation}\label{eq:reflections}
  \sigma_{X_1}\sigma_{X_2} = \sigma_{X_2}\sigma_{X_1} \iff X_1\perpgo X_2.
\end{equation}
One might call $\perpgo$ an orthogonality, but note that \eqref{eq:reflections} yields the formula
\begin{equation}\label{nibyizotrop}
  X_1\subset X_2 \implies  X_1 \perpgo X_2,
\end{equation} 
which fails to fit  intuitions 
that are commonly associated with the notion of an orthogonality of subspaces
in an Euclidean space.
For this reason we put some restrictions on $\perpgo$ to get a relation that 
conforms intuitions concerning Euclidean orthogonality more:
\begin{equation}\label{def:perpg}
  X_1\perpg X_2\;\; \defiff X_1\perpgo X_2 \text{ and } X_1\cap X_2 \neq X_1, X_2.
\end{equation}
So, in view of \eqref{nibyizotrop} the relations $\perpg$ and $\perpgo$ are closely related indeed:
\begin{equation}\label{g-go}
   X_1 \perpgo X_2 \iff X_1 \perpg X_2 \text{ or } X_1 \subset X_2 \text{ or }  X_2 \subset X_1. 
\end{equation}  
In view of \eqref{def:perpgo}, \eqref{wz:osieo} it is seen that the relation
$X_1 \perpg X_2$ can be 
characterized by any of the following three (mutually equivalent) conditions:
\begin{enumerate}\itemsep-2pt
\item[{(\ref{def:perpgo}.0)}]
  there is a point $q \in X_1\cap X_2$  and $Z_1, Z_2\in\hipy\setminus\hipy_0$ such that \\
  $q\in Z_1, Z_2 \perpx X_1 \cap X_2$, \,
  $Z_1 \perpx Z_2$, and  
  $(X_1 \cap X_2) \sqcup Z_i = X_i$ for $i=1,2$;
\item[{(\ref{wz:osieo}.$i$)}]
  there is  $Z_i\in\hipy\setminus\hipy_0$  such that \\ 
  $Z_i \perpx X_{3-i}$, \, $(X_1\cap X_2)\sqcup Z_i = X_i$, and not $X_{3-i}\subset X_i$;
\end{enumerate}
where $i=1,2$.

Let us write 
\begin{ctext}
  $X_1\perpu{m}{k_1}{k_2}X_2$
  when 
  $X_1\perpg X_2$,
  $X_1\in\hipy_{k_1}$,
  $X_2\in\hipy_{k_2}$, and 
  $X_1\cap X_2\in\hipy_{m}$.
\end{ctext}
Following this terminology we can say that the orthoadjacency relation $\oadjac_k$
considered in \cite{perpadja} is the relation $\perpu{k-1}{k}{k}$ for 
a fixed integer $k$.

Note the evident restrictions that dimensions $k_1,k_2,m$ must satisfy in order
to have $\perpu{m}{k_1}{k_2}$ nontrivial
\begin{equation}\label{zal:nontriv}
  \text{there are } X_1,X_2 \text{ such that }
  X_1 \perpu{m}{k_1}{k_2} X_2 
  \quad \iff \quad k_1 + k_2 - m \leq \dim({\goth M}).
\end{equation}

\begin{lem}\label{lem1}
  Let\/ $Y_1\in\hipy_{k_1-m}$ and\/ $X_2\in\hipy_{k_2}$ intersect in a point.
  Assume that $k_1 \le k_2$ and $k_1 + k_2 -m \le \dim({\goth M})$.
  The following conditions are equivalent
  \begin{enumerate}[\rm(i)]\itemsep-2pt
  \item\label{lem1:war1}
    $Y_1 \perpg X_2$ (i.e. actually, $Y_1 \perpx X_2$); 
  \item\label{lem1:war2}
    $X_1\perpg X_2$ for each $X_1\in\hipy_{k_1}$
    such that $Y_1\subset X_1$ and 
    $\dim(X_1\cap X_2)=m$.
  \end{enumerate}
\end{lem}
\begin{proof}
  The implication \eqref{lem1:war1}$\implies$\eqref{lem1:war2} follows directly from
  (\ref{wz:osieo}.1).
  
  Assume \eqref{lem1:war2}; set $V := Y_1\sqcup X_2$ and let $q\in Y_1\cap X_2$.
  Then $\dim(V) = k_1 + k_2 - m$.
  Let $W$ be the orthocomplement of $Y_1$ in $V$ through $q$, so 
  $\dim(W) = k_2$. Since $k_1 \le k_2$ we have 
  $m \le \dim(W\cap X_2)$, so there is $T\subset W\cap X_2$
  with $\dim({T}) = m$. Set $X_1: = T \sqcup Y_1$.
  Then $\dim(X_1) = k_1$ and thus $X_1\perpg X_2$. 
  Clearly, $X_1\cap X_2 = T$.
  By (\ref{wz:osieo}.1), there is $Z\in\hipy$ such that
  $Z\perpx X_2$ and $X_1 = T\sqcup Z$.
  Since both $Y_1, Z$ are orthocomplements of $T$
  in $X_1$, we get $Z = Y_1$ and thus \eqref{lem1:war1}
  follows.  
\end{proof}

\begin{lem}\label{lem2}
  Let\/ $1 \le k_1$ and\/ $1 < k_2$. Then for $L_1,L_2\in\lines$ the following
  conditions are equivalent
  \begin{enumerate}[\rm(i)]\itemsep-2pt
  \item
    $L_1\perp L_2$;
  \item
    there are $X_1\in\hipy_{k_1}$, $X_2\in\hipy_{k_2}$
    such that $X_1\perpx X_2$ and $L_i \subset X_i$
    for $i=1,2$.
  \end{enumerate}
\end{lem}

Notice that the assumption $1 < k_2$ in \ref{lem2} is significant as the lines
$L_1$, $L_2$ could be skew so, we need some more room in $X_2$ to find there
the translate of $L_2$ that meets $L_1$.

Now, let us consider the structure 
\begin{ctext}
  ${\goth K} := 
  \bstruct{S,\hipy_{k_1},\hipy_{k_2}, \mathord{\perpgo}\cap (\hipy_{k_1}\times\hipy_{k_2})}$;
\end{ctext}
for fixed $k_1,k_2$ such that $1\le k_1,k_2< \dim({\goth M})$.
As the inclusion relations involved in \eqref{g-go} and \eqref{def:perpg} are 
expressible in terms of pure incidence
language of \bstruct{S,\hipy_{k_1},\hipy_{k_2}} it is easily seen that
$\goth K$ and $\bstruct{S,\hipy_{k_1},\hipy_{k_2}, \mathord{\perpg}\cap (\hipy_{k_1}\times\hipy_{k_2})}$
are definitionally equivalent.
\begin{thm}\label{thm:main}
  Let $1 \le k_1,k_2$.
  \begin{sentences}\itemsep-2pt
  \item\label{main:1}
    If\/ $k_1 + k_2 - m \le \dim({\goth M})$, then
    the Euclidean space $\goth M$ is definable in the structure 
    $\struct{S,\hipy_{k_1},\hipy_{k_2},\perpu{m}{k_1}{k_2}}$.
  \item\label{main:2}
    The Euclidean space $\goth M$ is definable in $\goth K$.
  \end{sentences}
\end{thm}
\begin{proof}
  By \ref{fct:hipy<->lines}, for each integer $n$ the set 
  $\hipy_n$ is definable in 
  the reduct $\struct{S,\hipy_{k_1},\hipy_{k_2}}$ of $\goth K$.
  In particular, the family  $\lines$ of lines of $\goth M$ is definable in $\goth K$.
  Moreover, $\perpu{m}{k_1}{k_2}$ is definable in $\goth K$
  for each sensible $m$. Without loss of generality we can assume that
  $k_1 \le k_2$. By \ref{lem1}, the relation $\perpu{0}{k_1-m}{k_2}$
  is definable in $\goth K$ and in $\struct{S,\hipy_{k_1},\hipy_{k_2},\perpu{m}{k_1}{k_2}}$.
  Finally, by \ref{lem2} the proof is complete.
\end{proof}

\section{Synthetic properties of orthogonalities}

In this section we aim to show a few specific properties
of orthogonality relations $\perpgo$ and $\perpg$ considered on the family of all the subspaces of $\goth M$.
Some of them are analogous to known properties of the relation $\perp$
considered on the lines of $\goth M$, but there are also 
remarkable differences.

\subsection{Orthogonality $\perpg$}

\begin{fact}\label{fact:ax}
  Let\/ $A, B, C\in\hipy$. 
  \begin{sentences}
  \item\label{fact:ax:a}
    If\/ $A\perpg B$, then $B\perpg A$.
  \item\label{fact:ax:b}
    If\/ $A\perpg B$, then $A\cap B\neq\emptyset$.
  \item\label{fact:ax:c}
    If\/ $A\perpg B\parallel C$ and $A\cap C\neq\emptyset$, then 
    $A\perpg C$.
  \item\label{fact:ax:d}
    There are no nonempty $D_1, D_2\in\hipy\setminus\hipy_0$ with 
    $D_1\subseteq D_2$, and $D_1\perpg D_2$.
  \item\label{fact:ax::e}
    If\/ $\emptyset\neq A\subsetneq B\subsetneq C$, then
    there is the unique $B'\in\hipy$ such that 
    $B\cap B' = A$, $B\perpg B'$, and $B\sqcup B' = C$. 
  \end{sentences}
\end{fact}

\begin{prop}\label{prop:perpsqcup}
  Let\/ $A, B, C\in\hipy$. If\/ $A\perpg B$ and $A\perpg C$, then 
  $A\perpg(B\sqcup C)$ or $A\subseteq B\sqcup C$.
\end{prop}

\begin{proof}
  Assume that $A\perpg B$, $A\perpg C$, and $A\nsubseteq B\sqcup C$. 
  By \ref{fact:ax}\eqref{fact:ax:a} we have $A\cap B\neq\emptyset$.
  Since $A\cap B\subseteq A\cap(B\sqcup C)$ there is a common point
  $q$ of $A$ and $B\sqcup C$. From our assumption and (\ref{wz:osieo}.1)
  there are $Z_B, Z_C\in\hipy$ such that 
  \begin{equation}\label{eq:prop:perpsqcup:A}
    Z_B\perpx A,\quad (A\cap B)\sqcup Z_B = B\qquad\text{and}\qquad
      Z_C\perpx A,\quad (A\cap C)\sqcup Z_C = C.
  \end{equation}    
  Take $Z := Z_B'\sqcup Z_C'$, where $Z_B', Z_C'$ are translates of $Z_B, Z_C$ 
  respectively, through $q$. Therefore, by \eqref{eq:perpxsup} and 
  \eqref{eq:prop:perpsqcup:A} we have $A\perpx Z_B'\sqcup Z_C'=Z$.
  Now as $q\in A, B\sqcup C, Z$ and $Z\subseteq B\sqcup C$ we have 
  $Z\sqcup (A\cap (B\sqcup C)) = (Z\sqcup A)\cap (B\sqcup C)$. 
  Note that the equalities in \eqref{eq:prop:perpsqcup:A} give
  $Z_B'\sqcup A = Z_B\sqcup A = B\sqcup A$ and $Z_C'\sqcup A = Z_C\sqcup A = C\sqcup A$.
  So, we have $Z\sqcup A = (B\sqcup A)\sqcup (C\sqcup A) = (B\sqcup C)\sqcup A$
  and finally $Z\sqcup (A\cap (B\sqcup C)) = B\sqcup C$ which by (\ref{wz:osieo}.2)
  gives our claim.
\end{proof}

In some specific cases $\perpg$ may be transitive under inclusion which is showed
in next two propositions.
\begin{prop}\label{prop:cosik2}
  Let $A,B,C \in \hipy$. If\/ $A \perpg B$ and 
  $A\cap B \subsetneq C\subset B$, 
  then $A\perpg C$.
\end{prop}
\begin{proof}
  Let $q \in A\cap B$.  
  From assumptions, $A\cap C = A\cap B$.  
  By (\ref{wz:osieo}.1), there is $A'\in\hipy$ with
  $q \in A'$, $A = (A\cap B) \sqcup A'$ and $A'\perpx B$.
  Thus $A' \perpx C$ and the claim follows by (\ref{wz:osieo}.1).
\end{proof}
%
%

\begin{prop}\label{prop:cosik}
  Let\/ $A, B, C\in\hipy$. If\/ $A\perpg B$ and $A\cap B\subset C\subsetneq A$,
  then $A\perpg(B\sqcup C)$.
\end{prop}

\begin{proof}
  Let $q\in A\cap B$. Note that $A, B, C$ lay in the bundle through $q$, i.e.
  in a projective space, and thus we have $C = (A\cap B)\sqcup C = A\cap (B\sqcup C)$.
  In view of (\ref{wz:osieo}.1) 
  there is $Z\in\hipy$ such that $Z\perpx A$ and
  $B = (A\cap B)\sqcup Z$. From the latter equality we have
    $$B\sqcup C = (A\cap B)\sqcup Z\sqcup C = A\cap (B\sqcup C)\sqcup Z$$
  which, together with $Z\perpx A$, again by (\ref{wz:osieo}.1) completes the proof.
\end{proof}

The following example shows that it is hard to tell anything  more about transitivity 
of $\perpg$ than it is said in \ref{prop:cosik2} and \ref{prop:cosik}.

\begin{exm}\label{exm:notrans}\strut
  \begin{sentences}\itemsep-2pt
  \item\label{exm:notrans:1}
    There are $A,B,C\in \hipy$ such that
	\begin{equation*}
	  A \perpg B \subset C,\quad
          \neg A \perpg C, \quad\text{and}\quad 
          \dim(C) = \dim(B) + 1.
	\end{equation*}
  \item\label{exm:notrans:2}
    There are $A,B,C\in \hipy$ such that
	\begin{equation*}
	  A \perpg B \supset C ,\quad
          \neg A\perpg C,\quad
	  A\cap  C\neq\emptyset,
          \quad\text{and}\quad \dim(C) = \dim(B) - 1.
	\end{equation*}
  \end{sentences}
  In essence, one can take lines $A, B$ and a plane $C$ in \eqref{exm:notrans:1},
  as well as, planes $A,B$ and a line $C$ in \eqref{exm:notrans:2}.
\end{exm}

Therefore no ``simple'' form of transitivity can be proved. We finish with 
yet another property of $\perpg$.

\begin{prop}\label{prop:meet}
  Let\/ $A, B, C\in\hipy$. 
  If\/ $A\perpg B$, $A\perpg C$, and $A\cap B\cap C\neq\emptyset$, then 
  $A\perpg(B\cap C)$ or $B\cap C\subseteq A$.
\end{prop}
\begin{proof}
  We can assume that $B\cap C$ is at least a line as otherwise our claim 
  is clear. Let $q\in A\cap B\cap C$.
  Thanks to (\ref{wz:osieo}.1) we can take $Z_B, Z_C\in\hipy$ such that
  \begin{equation}\label{eq:prop:meet:A}
    Z_B\perpx B,\quad (A\cap B)\sqcup Z_B = A\qquad\text{and}\qquad
      Z_C\perpx C,\quad (A\cap C)\sqcup Z_C = A.
  \end{equation}   
  Note that $Z_B$ is the orthocomplement of $A\cap B$ in $A$ through $q$
  and $Z_C$ is the orthocomplement of $A\cap C$ in $A$ through $q$. So, 
  slightly abusing notation we can write
    $$Z := Z_B\sqcup Z_C = (A\cap B)^\perp\sqcup(A\cap C)^\perp = (A\cap B\cap C)^\perp.$$
  Hence $Z\sqcup (A\cap B\cap C) = A$. Moreover 
  $q\in Z_B, Z_C\perpx B\cap C\ni q$ by \eqref{eq:prop:meet:A}. Hence 
  by \eqref{eq:perpxsup} we get $Z\perpx B\cap C$, which in view of (\ref{wz:osieo}.1)
  suffices as a final argument.
\end{proof}

\subsection{Orthogonality $\perpgo$}

According to \eqref{def:perpg} or \eqref{g-go}, properties of relation $\perpgo$ are
simple consequences of properties of relation $\perpg$ with possible inclusions between
its arguments taken into account.

\begin{prop}\label{fact:axo}
  Let\/ $A, B, C\in\hipy$. 
  \begin{sentences}
  \item\label{fact:ax:ao}
    If\/ $A\perpgo B$, then $B\perpgo A$.
  \item\label{fact:ax:bo}
    If\/ $A\perpgo B$, then $A\cap B\neq\emptyset$.
  \item\label{fact:ax:co}
    If\/ $A\perpgo B\parallel C$ and $A\cap C\neq\emptyset$, then 
    $A\perpgo C$.
  \item\label{fact:ax:do}
    If\/ $\emptyset\neq A\subsetneq B\subsetneq C$, then
    there is the unique $B'\in\hipy$ such that 
    $B\cap B' = A$, $B\perpgo B'$, and $B\sqcup B' = C$. 
  \item\label{fact:ax:eo}
    If\/ $A\perpgo B$ and $A\perpgo C$, then 
    $A\perpgo(B\sqcup C)$.   
  \item\label{fact:ax:fo}      
    If\/ $A \perpgo B$ and 
    $A\cap B \subset C\subset B$, then $A\perpgo C$.
  \item\label{fact:ax:go}
    If\/ $A\perpgo B$ and $A\cap B\subset C\subset A$,
    then $A\perpgo(B\sqcup C)$.    
  \item\label{fact:ax:ho}
    If\/ $A\perpgo B$, $A\perpgo C$, and $A\cap B\cap C\neq\emptyset$, then 
    $A\perpgo(B\cap C)$.
  \end{sentences}
\end{prop}
\begin{proof}
  \eqref{fact:ax:ao} -- \eqref{fact:ax:do} 
  follow directly from \ref{fact:ax} and \eqref{g-go}.
  
  \eqref{fact:ax:eo}: 
  It suffices to apply \eqref{g-go} plus \ref{prop:perpsqcup} or \eqref{nibyizotrop}.
  Only two cases: (a) $C \subset A \perpg B$, (b) $B \subset A \perpg C$ 
  of interpretation of the assumptions may appear problematic,
  but they are equivalent up to names of variables.
  Assume that (a) holds. Set $C' := C \sqcup (A \cap B)$. Then  
  $C\sqcup B = C\sqcup (A\cap B)\sqcup B = C'\sqcup B$.
  So, we have $A \cap B \subset C' \subset A$.
  If $C' = A$, then the conclusion of \eqref{fact:ax:eo} follows by \eqref{nibyizotrop}.
  If $C' \neq A$ the claim follows by \eqref{prop:cosik}.
    
  \eqref{fact:ax:fo} 
  is immediate by \eqref{g-go} plus \ref{prop:cosik2} or \eqref{nibyizotrop}.
  
  \eqref{fact:ax:go} 
  is immediate by \eqref{g-go} plus \ref{prop:cosik} or \eqref{nibyizotrop}.
  
  \eqref{fact:ax:ho}: 
  Apply \eqref{g-go} plus \ref{prop:meet} or \eqref{nibyizotrop}.
  Two cases, though equivalent up to variables, of the assumptions may raise some problems:
  (a) $B\perpg A\subset C$ (b) $C\perpg A\subset B$. Assume that (a) holds.
  Set  $C' := B\cap C$. Then $A\cap B\subset C' \subset B$.
  If $A \cap B \neq C'$, then the claim comes from \ref{prop:cosik2}.
  If $A \cap B = C'$, then $B\cap C = C' = A \cap B\subset A$ and the claim is a consequence
  of \eqref{nibyizotrop}.
\end{proof}

\medskip
\par\noindent
Author's address:
\par\noindent
Jacek Konarzewski, Mariusz {\.Z}ynel
\\
Institute of Mathematics, University of Bia{\l}ystok
\\
Akademicka 2, 15-267 Bia{\l}ystok, Poland
\\
e-mail: \verb+konarzewski20@wp.pl+, \verb+mariusz@math.uwb.edu.pl+

\end{document}